\documentclass[a4paper, 12pt]{article}
\usepackage{calrsfs}
\usepackage{tikz-cd}
\usepackage[utf8x]{inputenc}
\usepackage[T1]{fontenc}
\usepackage[top=2cm, bottom=2cm, left=2cm, right=2cm]{geometry}
\usepackage{amssymb}
\usepackage{amsmath}
\usepackage[backend=biber]{biblatex}
\newtheorem{lem}{Lemma}
\newtheorem{prop}{Proposition}
\newtheorem{thm}{Theorem}
\newtheorem{cor}{Corollary}
\newtheorem{defi}{Definition}

\author{Antonin Assoun}
\title{Regular version of the inverse Galois problem for skew fields}
\nocite{*}
\begin{document}
\maketitle
\begin{abstract}
In order to extend the study of the regular version of the inverse Galois problem to skew fields, we generalize the definition of regular field extensions for commutative fields to the case of arbitrary fields. We then propose a general version of the regular inverse Galois property and show that for a field $k$, the study of the class $\mathrm{InvReg}(k)$ of non-trivial regular extensions of $k$ satisfying this property constitutes a natural generalization of the classical regular inverse Galois problem. Next, we use recent results from Behajaina, Deschamps, and Legrand on the inverse Galois problem in the noncommutative setting to show that certain skew fraction fields $h(t,\sigma)$ belong to the class $\mathrm{InvReg}(h)$. Finally, we generalize Behajaina’s method to construct extensions belonging to $\mathrm{InvReg}(h)$ that are not of the form $h(t,\sigma)$ with $\sigma$ an automorphism of finite order.
\end{abstract}

\section{Introduction}
The Galois theory of commutative fields admits a generalization to arbitrary fields, developed by Jacobson (see \cite{Jac}, \cite{Cohn_1995}). This allows in particular the formulation of the inverse Galois problem for arbitrary fields. Recent work on the subject by Behajaina, Deschamps, and Legrand (\cite{BEH}, \cite{Des21}, \cite{Des24}, \cite{DesLe}) leads to the consideration of the problem 
\newline
\newline
$\mathbf{IGP}_k$ \textbf{(Inverse Galois Problem over an arbitrary field $k$)} : \textit{Given a field $k$ we say that the property} $\mathrm{IGP}_k$\ \textit{is satisfied if every finite group is the Galois group of an outer\footnote{An extension $h/k$ is outer if the only inner automorphism of $h$ restricting to the identity on $k$ is the identity on $h$. These extensions are those for which the Galois correspondence is identical to the correspondence for commutative fields.} Galois extension of $k$.}
\
\newline
\newline
This is a generalization of the classical inverse Galois property for commutative fields. The goal of this article is to study a generalization of the regular variant of this problem:
\newline
\newline
$\mathbf{RIGP_k}$ \textbf{(Regular Inverse Galois Problem over a commutative field $k$)}: \textit{For every finite group} $G$, \textit{there exists a commutative Galois extension} $L/k(t)$ \textit{regular over} $k$ \textit{with} $\mathrm{Gal}(L/k(t)) \simeq G$.
\newline
\newline
The study of this variant has led to profound results in field arithmetic. For example, Pop's work \cite{Pop96} shows that this property holds for ample fields\footnote{An ample field is a commutative field $k$ such that every absolutely irreducible curve over $k$ having a smooth rational point has infinitely many rational points. This property is equivalent to the fact that $k$ is existentially closed in its field of Laurent series. Henselian fields and algebraically closed fields are ample.}. Ample fields even satisfy the stronger property: if $k$ is an ample field and $G$ is a finite group, there exists a field $L$ with $k(t) \subset L \subset k((t))$ such that $L/k(t)$ is a Galois extension with group $G$. Following Behajaina, we will denote this variant by $\mathrm{SIGP}_k$ with an S for series.
\newline
\newline
In section 2, we propose a generalization to arbitrary fields of the notion of regular extension of a field (definition 2). If $K/k$ is a regular extension in this sense, we then define the regular inverse Galois problem relative to $K/k$:
\begin{defi}[$\mathbf{RIGP}_{K/k}$ Regular Inverse Galois Problem for a regular extension $K/k$]
We say that the property $\mathrm{RIGP}_{K/k}$ holds if for every finite group $G$, there exists an outer Galois extension $L/K$, regular over $k$, whose Galois group is isomorphic to $G$.
\end{defi}
We then denote by $\mathrm{InvReg}(k)$ the class of regular extensions $K/k$ for which $\mathrm{RIGP}_{K/k}$ holds. When $k$ is a commutative field, it suffices that the extension $k(t)/k$ belongs to the class $\mathrm{InvReg}(k)$ to guarantee that every regular commutative extension of $k$ of finite type belongs to $\mathrm{InvReg}(k)$ (proposition 3)\footnote{The hypothesis of proposition 3 is more general but is true for finite type extensions.}. This reduces the regular inverse Galois problem for commutative fields $k$ to a problem on the field $k(t)$. Such a reduction to the extension $k(t)/k$ is unclear in the non-commutative context: whether a regular extension $K/k$ of finite type in the class ${\rm InvReg}(k)$ may depend on $K$ (and not 
just on $k$). 
\newline
In section 2, we investigate whether the skew fraction fields $k(t,\sigma)$ of rational fractions of polynomials twisted by an automorphism $\sigma$ of $k$ belong to the class $\mathrm{InvReg}(k)$. We refer to \cite{Cohn_1995} or \cite{Ore} for the definitions of these skew fraction fields and their analogs $k((t,\sigma))$, the twisted Laurent series fields. We show that recent work by Behajaina, Deschamps, and Legrand (\cite{BEH}, \cite{Des21},\cite{DesLe}) about realizing finite groups as Galois groups of extensions of $k(t,\sigma)$ (or $k(t,\sigma,\delta)$ in \cite{Des21}) implies that $k(t,\sigma) \in \mathrm{InvReg}(k)$, i.e., $\mathrm{RIGP}_{k(t,\sigma)/k}$ holds, when $k$ is a field which center contains a subfield fixed by $\sigma$ and ample.

In sections 3 and 4, we provide new examples of extensions $K/k$ for which the property $\mathrm{RIGP}_{K/k}$ holds. To do this, we generalize the scalar extension method used in \cite{BEH} and \cite{Des21} to apply it to cases where $\sigma$ is not necessarily an automorphism of finite order but only \textit{locally} of finite order: that is, for every $x \in K$, there exists a non-zero natural number $n$ such that $\sigma^n(x) = x$. This allows us to prove the following result:
\begin{thm}
Let $h/k$ be an algebraic outer Galois extension of fields such that $\mathrm{Gal}(h/k)$ is an infinite product of cyclic groups of coprime orders. If $[Z(h):Z(k)] = +\infty$, and $Z(k)$ contains an ample field, then the class $\mathrm{InvReg}(h)$ contains extensions of $h$ that are non-commutative and not of the form $h(t,\sigma)/h$ with $\sigma$ an automorphism of finite order.
\end{thm}
For an example of an extension $h/k$ satisfying the hypotheses of theorem 1, let $\Gamma$ be the direct product of all cyclic groups of prime order, let $k$ be an ample field and let $L/k(t)$ be a commutative Galois extension of group $\Gamma$. Then $h=L$ and $k=k(t)$ satisfy the hypotheses. This theorem is obtained by considering inductive limits of extensions constructed in \cite{BEH}: we consider a filtration of $h$ by subfields $h_n$ stable under a topological generator $\sigma$ of $\mathrm{Gal}(h/k)$, such that $\sigma$ is of finite order on each $k_n$, and show that an inductive limit of extensions of $h_n$ obtained using the method in \cite{BEH} gives a suitable extension of $h$. An important tool will be a generalization to arbitrary fields of Artin's lemma for infinite algebraic extensions as stated in \cite{Fr}. 
\newline
\newline
\textbf{Notation:} From now on, we adopt the convention that a field $\mathbf{k}$ written in bold is commutative, that if two fields $K,k$ are written with the same letter, the one in uppercase is a regular extension of the one in lowercase, and for a ring $A$, we denote by $Z(A)$ its center.
\section{Regularity of Field Extensions and Regular Inverse Galois Problem}
In this section, we extend the property of regularity of a commutative field extension to skew fields. We explain the differences with the commutative case and give some examples of regular extensions of arbitrary fields.
\subsection{Regular Extensions}
In this subsection, we give the definition of a regular extension that we use later and discuss some properties of these extensions.
\newline
An element $x$ of a field extension $h/k$ is said to be algebraic over $k$ if the subfield it generates, $k(x) \subset h$, has finite dimension as a left \textit{and} right vector space over $k$\footnote{This definition is the most convenient to use for Galois theory. However there are alternative generalizations of algebraicity which are not known to be equivalent. One could ask for the property that any finite family of elements of $h$ generates a finite extension of $k$. Or one could ask that for any $x$ in $h$ there exists a polynomial (or even left or right polynomial) equation of which $x$ is a zero.}. Although it is possible to have extensions which have distinct left and right degrees, for Galois extensions, the left and right degrees are always equal (see \cite{Cohn_1995}, corollary 3.3.8). In particular, if $h/k$ is a field extension, $x \in h$, and $k(x)$ is contained in a sub-extension $l/k$ that is Galois of finite degree, then $x$ is algebraic over $k$. The following formulation of the definition of a regular extension is valid for any field:
\begin{defi}
We say that a field extension $K/k$ is regular or that $K$ is a regular extension of $k$ if every element of $K$ algebraic over $k$ is an element of $k$.
\end{defi}
When $K$ and $k$ are commutative, it is equivalent to say that if $\overline{K}$ is an algebraic closure of $K$ and $\overline{k}$ is the algebraic closure of $k$ in $\overline{K}$, then $K \cap \overline{k} = k$. This definition does not extend naturally: first, a skew field $K$ may have several non-isomorphic maximal algebraic extensions (see \cite{Leg}). Furthermore, if $\overline{K}$ is \textit{one} such maximal algebraic extension of $K$, it is not known whether the set of elements of $\overline{K}$ algebraic over $k$ forms a field.
\newline
\newline
Let us mention a property of regular extensions that will be an important tool to prove the regularity of an extension:
\begin{prop}
Let $k \subset L \subset K$ be fields. If $K/k$ is a regular extension, then $L/k$ is a regular extension.
\end{prop}
\textbf{Proof:} If $x \in L$ is algebraic over $k$, then since $k(x) \subset K$, $x$ is an element of $K$ algebraic over $k$. Since $K$ is a regular extension of $k$, we have $x \in k.\blacksquare$
\newline
\newline
However, it is not clear that if $K/k$ is regular and $H/K$ is regular, then $H/k$ is regular. When $H$ is a commutative field, this property holds, but its proof relies on the fact that the evaluation map in $x \in H$ of a polynomial $k[X] \to H$ is a ring homomorphism, which is only true for $x \in Z(H)$. When $k$ is an arbitrary field, $A$ is a $k$-ring, and $x \in A$, if $P(t) = \sum_{i=0}^n a_n t^i \in k[t]$, we will conventionally write $P(x) = \sum_{i=0}^n a_i x^i \in A$ as the right evaluation of $P$ at $x$. This evaluation satisfies the following property: if $f,g \in k[t]$, if $y \in A$ and $g(y) \neq 0$, then:
$$ fg(y) = f(g(y)xg(y)^{-1})g(y). $$ 
Here, for $x \in H$, the existence of a polynomial $P \in k[t]$ such that $P(x) = 0$ remains a necessary condition for $x$ to be algebraic over $k$, but it is no longer a sufficient condition, and the proof used in the commutative case does not extend to this setting.
\subsection{Existentially Closed Extensions}
We extend to fields of finite dimension over their center a property known for commutative fields: if $L/K$ is a field extension such that $K$ is existentially closed in $L$, then $L/K$ is a regular extension. The result is the following.
\begin{thm}
Let $K$ be a field of finite dimension over its center. If $L/K$ is a field extension such that $K$ is existentially closed in $L$ and $Z(K) \subset Z(L)$, then $L/K$ is a regular extension.
\end{thm}
To prove the theorem we need the following lemma.
\begin{lem}
Let $K$ be a field of finite dimension over its center, and $L/K$ a finite degree field extension such that $Z(K) \subset Z(L)$, $x$ an element of $L$, $N \geq 1$ a natural number, and $a_0, \ldots , a_{n-1}$ elements of $K$. Then the formula:
$$ \exists y \in L \ xy = yx \wedge y^N + \sum_{i=0}^{N-1} a_i y^i = 0 $$
is equivalent to an existential sentence in the language of $K$-rings.
\end{lem}
\textbf{Proof:} Choose a basis $(e_i)_{1 \leq i \leq n}$ of $L$ as a $Z(K)$-vector space. There exist elements $(\lambda^{ij}_k)_{1 \leq i,j,k \leq n}$ in $Z(K)$ such that for $1 \leq i,j \leq n$ we have $$e_i e_j = \sum_{k=1}^n \lambda^{ij}_k e_k.$$ 
Also, there exist elements $x_1, \ldots, x_n \in Z(K)$ such that $$x = \sum_{i=1}^n x_i e_i.$$
Thus, if $y = \sum_{j=1}^n y_j e_j$, then:
\begin{align*}
xy = yx &\Longleftrightarrow \sum_{i,j=1}^n x_i y_j e_i e_j = \sum_{i,j=1}^n x_i y_j e_j e_i\\
&\Longleftrightarrow \sum_{i,j,k=1}^n x_i y_j \lambda_k^{ij} e_k = \sum_{i,j,k=1}^n x_i y_j \lambda_k^{ji} e_k.
\end{align*}
For $k \in \mathbb{N}^*$, let $S_{k,n}$ denote the set of maps $\sigma : \{1, \ldots, k\} \to \{1, \ldots, n\}$. For $\sigma \in S_{k,n}$, there exist constants $\lambda_j^{\sigma}$ in $Z(K)$ such that $e_{\sigma(1)} \ldots e_{\sigma(n)} = \sum_{j=1}^n \lambda_j^{\sigma} e_j.$ We then have:
$$ y^k = \sum_{j=1}^n \left( \sum_{\sigma \in S_{k,n}} \lambda_j^{\sigma} y_{\sigma(1)} \ldots y_{\sigma(k)} \right) e_j.$$
By breaking down the $a_k$ for $0 \leq k \leq n$ on the basis $(e_j)_{1 \leq j \leq n}$, we obtain $a_k = \sum_{l=1}^n a_l^k e_l$, and hence:
\begin{align*}
y^N + \sum_{k=0}^{N-1} a_k y^k &= \sum_{j=1}^n \left( \sum_{\sigma \in S_{N,n}} \lambda_j^\sigma y_{\sigma(1)} \ldots y_{\sigma(N)} \right) e_j \\
&+ \sum_{k=1}^{N-1} \left( \sum_{l=1}^n a_l^k  \left( \sum_{j=1}^n \left( \sum_{\sigma \in S_{k,n}} \lambda_j^\sigma y_{\sigma(1)} \ldots y_{\sigma(k)} \right) \sum_{t=1}^n \lambda_t^{lj} e_t \right)\right)\\
&+ \sum_{l=0}^n a_l^0 e_l \\
&= \sum_{t=1}^n \left( a_t^0 + \sum_{\sigma \in S_{N,n}} \lambda_t^\sigma y_{\sigma(1)} \ldots y_{\sigma(N)} \right) e_t + \sum_{t=1}^n\left(\sum_{k=1}^{N-1} \sum_{l=1}^n  \sum_{\sigma \in S_{k,n}} a_l^k \lambda_t^\sigma \lambda_j^{lt} y_{\sigma(1)} \ldots y_{\sigma(k)}\right) e_t.\\
&=\sum_{t=1}^n \left( a_t^0 + \sum_{\sigma \in S_{N,n}} \lambda_t^\sigma y_{\sigma(1)} \ldots y_{\sigma(N)}+\sum_{k=1}^{N-1} \sum_{l=1}^n  \sum_{\sigma \in S_{k,n}} a_l^k \lambda_t^\sigma \lambda_j^{lt} y_{\sigma(1)} \ldots y_{\sigma(k)}\right) e_t
\end{align*}
Each of the formulas $$\sum_{i,j=1}^n (\lambda_k^{ij} -\lambda_k^{ji})x_i y_j=0$$ for $1\leq k\leq n$ and $$a_t^0 + \sum_{\sigma \in S_{N,n}} \lambda_t^\sigma y_{\sigma(1)} \ldots y_{\sigma(N)}+\sum_{k=1}^{N-1} \sum_{l=1}^n  \sum_{\sigma \in S_{k,n}} a_l^k \lambda_t^\sigma \lambda_j^{lt} y_{\sigma(1)} \ldots y_{\sigma(k)}=0$$ for $1\leq t\leq n$ are quantifier free formulas in the language of $K$-rings with free variables $y_1,\ldots, y_n$. Therefore, so is their conjunction $\mathfrak{P}$ and $$\exists y_1\ldots\exists y_n \mathfrak{P}$$ is an existential sentence in the language of $K$-rings, equivalent to $$ \exists y \in L \ xy = yx \wedge y^N + \sum_{i=0}^{N-1} a_i y^i = 0 .\blacksquare$$
\textbf{Proof of theorem 2:} Let $x \in L$ be an element algebraic over $K$. Then there exists a unique monic polynomial of minimal degree $\pi_x = X^N + \sum_{i=0}^{N-1} a_i X^i$ such that $x$ is a right root of this polynomial, i.e., $x^N + \sum_{i=0}^{N-1} a_i x^i = 0$. By Lemma 1, the following formula is equivalent to an existential formula with coefficients in $K$:
$$ \exists y \, y^n + \sum_{i=0}^{n-1} a_i y^i = 0 \wedge xy = yx. $$
Now, since $x$ satisfies this formula, and as $K$ is existentially closed in $L$, there exists an element $x_0 \in K$ that satisfies this formula. In particular, $x_0$ is a right root of $\pi_x$, and thus there exists a polynomial $p(t)$ in the ring of central indeterminate polynomials $K[t]$ such that $\pi_x(t) = p(t)(t - x_0) \in K[t]$.
\newline
In particular, since $\pi_x(x) = 0$, we either have $(t - x_0)(x) = 0$, meaning $x = x_0$, or $\pi_x(x) = p((x - x_0)x(x - x_0)^{-1})(x - x_0) = 0$. Therefore, we must have $p((x - x_0)x(x - x_0)^{-1}) = 0$, but since $x$ and $x_0$ commute, this means that $p(x) = 0$. However, $p(t)$ has degree strictly smaller than $\pi_x$, which contradicts the minimality of the degree of $\pi_x$. Thus, $x = x_0$, i.e., $x \in K$.$\blacksquare$
\subsection{Skew Series Fields}
We show that skew series fields are regular extensions of their coefficient fields and provide an example of transcendental extensions of a field $k$ that only have $k$ as a common subextension up to $k$-isomorphism.
\newline
Let us recall the definitions of the skew series fields we use. Given a field $k$, an automorphism $\sigma$ of $k$, and a $\sigma$-derivation $\delta$ of $k$, we define a skew series field with right-hand coefficients $k_{\sigma,\delta}((z))$ (see \cite{Cohn_1995} section 2.3 for more details on the construction of the Cauchy product from the commutation relation and the properties of these fields). Its elements are formal series $S$ with coefficients in $k$ such that there exists an integer $n \in \mathbb{Z}$ with $S = \sum_{i=n}^{+\infty} z^i s_i$, and the multiplication on the left by an element $a$ of $K$ is given by the relation $az = \sum_{i=1}^{+\infty} z^i \sigma \circ \delta^{i-1}(a)$.
\newline
The skew series field denoted $k((t,\sigma))$ consists of formal series $S$ with coefficients in $k$ such that there exists an integer $n \in \mathbb{Z}$ with $S = \sum_{i=n}^{+\infty} s_i z^i$, and the multiplication on the right by an element $a$ of $k$ is given by the relation $z a = \sigma(a) z$.
\newline
We then have $k$-embeddings $k(t,\sigma,\delta) \to k_{\sigma,\delta}((z))$ given by $t \mapsto z^{-1}$ and $k(t,\sigma) \to k((t,\sigma))$ given by $t \mapsto t$.
\newline
These fields are generalizations of Laurent series fields: if $k$ is commutative, $\sigma = \mathrm{Id}_k$, and $\delta = 0$, then we recover $k_{\sigma,\delta}((z)) \simeq k((t,\sigma)) \simeq k((t))$.

\begin{prop}
    Let $k$ be a field, $\sigma$ an automorphism of $k$, and $\delta$ a $\sigma$-derivation of $k$. Then:
    \begin{itemize}
        \item The skew series field $k((t,\sigma))$ is a regular extension of $k$.
        \item The skew series field $k_{\sigma,\delta}((z))$ is a regular extension of $k$.
        \item The skew rational fraction field $k(t,\sigma,\delta)$ is a regular extension of $k$. 
    \end{itemize}
\end{prop}
\textbf{Proof:} Let $S(z) \in k_{\sigma,\delta}((z))$. Then there exist elements $(s_n)_{n \in \mathbb{Z}}$ in $k$ such that $S(z) = \sum_{-\infty}^{\infty} z^n s_n$. If $S(z) \notin k$, there exists an element $S'(z) \in k_{\sigma,\delta}((z))$ with non-zero $z$-adic valuation such that $k(S(z)) = k(S'(z))$. Indeed, if $S(z)$ has non-zero $z$-adic valuation, it suffices to choose $S'(z) = S(z)$; otherwise, $S(z) - s_0$ is non-zero and has non-zero $z$-adic valuation, so we can choose $S'(z) = S(z) - s_0$. In particular, if $n$ and $m$ are distinct integers, $S'(z)^n$ and $S'(z)^m$ have distinct $z$-adic valuations, and the family $(S'(z)^n)_{n \in \mathbb{Z}}$ is $k$-linearly independent (both left and right). Therefore, $S(z)$ is not algebraic. Replacing $z$ with $t$ in the previous proof gives a proof for $k((t,\sigma))$. Since there exists a $k$-embedding $k(t,\sigma,\delta) \to k_{\sigma,\delta}((z))$, it follows from Proposition 1 that $k(t,\sigma,\delta)$ is a regular extension of $k$.$\blacksquare$
\newline
\newline
\textbf{Remark:} This result also extends to more general series fields. The proof uses the properties of $\mathbb{Z}$ endowed with the usual order as a totally ordered group. For any well-ordered subset $S$ of a totally ordered group $G$ and any natural integer $n \in \mathbb{N}^*$, we have that $\mathrm{min}(S^n) = \mathrm{min}(S)^n$. The essential property of $\mathbb{Z}$ for the proof of Proposition 2 is that the only expression of $\mathrm{min}(S)^n$ as a product of $n$ elements of $S$ is $\mathrm{min}(S)^n$. Thus, our proof generalizes to any Malcev-Neumann series field (introduced by Neumann in \cite{Neu}, see also \cite{Cohn_1995} section 2.4) indexed by a totally ordered group $G$ such that for every element $g \in G$, the following property holds: for every $n \in \mathbb{N}$, if $g_1 \dots g_n = g^n$ and $g_1, \dots, g_n \geq g$, then $g_1 = \dots = g_n = g$. For example, for a field $k$, the field $k(\mathbb{Q})$ of Puiseux series over $k$ with central indeterminate is a regular extension of $k$.
\subsection{Regular Form of the Inverse Galois Problem}
Any regular extension $K/k$ with $K$ commutative has a subextension isomorphic to $k(t)$, and if $L/k(t)$ is a regular Galois extension with Galois group $G$, then $K \otimes_{k(t)} L/K$ is a regular Galois extension of $K$ with Galois group $G$. In particular, if every finite group occurs as the Galois group of a regular extension of $k(t)$, it also occurs as the Galois group of a regular extension of $K$ for any regular extension $K/k$. When the fields involved are non-commutative, there generally does not exist a regular extension $H/k$ such that every regular extension $L/k$ contains a subextension $k$-isomorphic to $H$. For example, if $k$ is a commutative field with an automorphism $\sigma$ of infinite order, then Lüroth's theorem asserts that every subextension of $k(t)/k$ is isomorphic to $k(t)$. However, no subextension of $k(t,\sigma)/k$ is isomorphic to $k(t)$. Indeed, if such a subextension existed, there would be an element of the center of $k(t,\sigma)$ not in $k$. In particular, using the $k$-embedding $k(t,\sigma) \to k((t,\sigma))$, we would have a series $S(t) \in k((t,\sigma))$ that commutes with all elements of $k$. But if $S(t) = \sum_{n=-\infty}^{\infty} a_n t^n$, using the commutativity of $k$, we see that saying $S(t)$ commutes with all elements of $k$ is equivalent to saying that for all $x \in k$ and for all $n \in \mathbb{Z}$ with $a_n \neq 0$, we have $\sigma^n(x) = x$. Since $\sigma$ has infinite order, the only non-zero coefficient of $S(t)$ can be $a_0$, and thus $S(t) \in k$. We have thus given an example of two non-trivial regular extensions $K_1/k$ and $K_2/k$ such that the only pair $(L_1, L_2)$ with $L_1/k$ a subextension of $K_1/k$, $L_2/k$ a subextension of $K_2/k$, and $L_1/k$ $k$-isomorphic to $L_2/k$ is the pair $(k, k)$. Said concisely, given regular extensions $K_1/k$ and $K_2/k$, the following diagram of extensions does always exist when $L_1$ and $L_2$ are commutative (and $K$ can be chosen $k$-isomorphic to $k(t)$) and it does not necessarily exist when they are skew fields .

$$\begin{tikzcd}
K_1 &  &                                &  & K_2 \\
    &  &                                &  &     \\
    &  & K \arrow[lluu] \arrow[rruu] &  &     \\
    &  &                                &  &     \\
    &  & k \arrow[uu]                   &  &    
\end{tikzcd}$$

Thus, a formulation of the regular inverse Galois problem over $k$ cannot reduce to the realization of finite groups as Galois groups of a unique regular extension of $k$, which is why we have defined in the introduction, for $K/k$ a field extension, the property $\mathrm{RIGP}_{K/k}$. In this paragraph, we study the relationship between this property and the property $\mathrm{RIGP}_k$ for $k$ a commutative field. Let us introduce classes of field extensions:
\begin{defi}
    Let $k$ be a field.
    \begin{itemize}
        \item We denote by $\mathrm{Reg}(k)$ the class of regular extensions of $h$ distinct from $h$.
        \item We denote by $\mathrm{InvReg}(k)$ the class of regular extensions $K/k$ such that $\mathrm{RIGP}_{K/k}$ is satisfied.
    \end{itemize}
\end{defi}
The following proposition summarizes relations between the property $\mathrm{RIGP}_k$ for $k$ commutative, the property $\mathrm{RIGP}_{K/k}$ for a regular field extension, and the classes introduced in the previous definition. It shows in particular that the study of the class $\mathrm{InvReg}(k)$ for $k$ an arbitrary field is a generalization to arbitrary fields of the study of the property $\mathrm{RIGP}_k$ when $k$ is commutative.
\newpage
\begin{prop}
If $k$ is a commutative field, then we have:
\newline
i) $\mathrm{RIGP}_{k(t)/k} \Longleftrightarrow \mathrm{RIGP}_k$
\newline
    ii) $\mathrm{RIGP}_k $ implies any regular extension $E/k$ of $k$, of finite degree over a purely transcendental extension of $k$ and monogenous\footnote{For commutative fields, regular extensions are usually required to be separable. In which case this hypothesis would be superfluous. But no fully convincing notion of separability for extensions of skew fields has been defined yet. We therefore chose to remove this hypothesis and to consider only monogenous extensions (which definition does not depend on the commutative setting) when needed. } over this extension belongs to $\mathrm{InvReg}(k)$.
\end{prop}
\textbf{Proof:}
i) If $L/k(t)$ is an outer Galois extension with finite Galois group, then $[L:k(t)] < +\infty$. Since the extension is outer, $k(t) \subset Z(L)$. Consequently, $L$ is a central simple algebra over its center. In particular, $L/Z(L)$ is a finite-degree Galois extension, and its Galois group is a subgroup of $\mathrm{Gal}(L/k(t))$. Now, by the Skolem-Noether theorem, we know that $\mathrm{Gal}(L/Z(L))$ consists of inner automorphisms of $L$. Since it is a subgroup of the Galois group of the outer extension $L/k(t)$, it must be the trivial group, and the Galois correspondence then gives $L = Z(L)$, meaning that $L$ is a commutative field. Therefore, if every finite group $G$ is isomorphic to the Galois group of an outer regular Galois extension, then every finite group $G$ is isomorphic to the Galois group of a regular Galois extension of commutative fields, and the converse is true because every commutative field extension is outer, which proves this point.
\newline
ii) Let $E/k$ be an extension as in the hypotheses of the proposition. Let $(x_i)_{i\in I}$ be a transcendence base of $E/k$ and $K=k((x_i)_{i\in\mathbb{I}})$. It is well known that $\mathrm{RIGP}_k$ implies that for any finite group $G$ there exists an extension $L/K$ which is Galois of group $G$ and regular over $k$. Applying this property for $G^n$ for $n\in\mathbb{N}^*$, one gets there are infinitely many pairwise linearly disjoint such $G$-Galois extensions. Let $(L_j/K)_{j\in\mathbb{N}}$ be a family of such extensions. Then for any $j\in\mathbb{N}$, $L_j\cap E$ is a subextension of $E/K$. As $L_j/K$ is Galois, $L_j$ and $E$ are linearly disjoint over $K$ if and only if their intersection is $K$. But as the $L_j$ are pairwise linearly disjoint over $K$, the subextensions $L_j\cap E$ (consider $E$ and the $L_j$ embedded in an algebraic closure $\Omega$ of $K$ so that the intersection makes sense) which are not trivial are all distinct. As $E/K$ is monogenous, it has only finitely many subextensions. Therefore infinitely many of the $L_j$ are linearly disjoint from $E$. Therefore there exists a $G$-Galois extension $L/K$, regular over $K$ and linearly disjoint from $E$. Then $E\otimes_K L$ is a field, and $E\otimes_KL/E$ is a Galois extension with Galois group $G$ which is a regular extension of $k$. Therefore $E$ belongs to $\mathrm{InvReg(k)}$.
\newline
\newline
\textbf{Remark :} For $ii)$ it would not be sufficient to assume $E/k$ is any regular extension .For instance if $k=\mathbb{C}$, $\mathrm{RIGP}_k$ is true.  But if $E$ is the field of Puiseux series over $\mathbb{C}$, it is a regular extension of $\mathbb{C}$ which cannot satisfy $\mathrm{RIGP}_{E/k}$ because it is algebraically closed and has no nontrivial commutative extension of finite degree.
\newline
\newline
When $k$ is a commutative field, we know that if it is Hilbertian, then $\mathrm{RIGP}_k \Rightarrow \mathrm{IGP}_k$. It is clear that this implication is not always satisfied: $p$-adic fields are ample and therefore satisfy the regular inverse Galois property, but all Galois groups over a $p$-adic field are solvable. Similarly, algebraically closed fields are ample, but only the trivial group can be realized as a Galois group over an algebraically closed field. Currently, we do not know of any generalization of Hilbert's property that applies to all fields, but it would be of special interest to know if such a property exists and allows the specialization of regular exterior Galois extensions.$\blacksquare$

\subsection{Interpretation of the work of Behajaina, Deschamps, and Legrand}

In this subsection, we show that the articles \cite{BEH}, \cite{Des21} provide examples of skew fields $k$, automorphisms $\sigma$ of $k$, and $\sigma$-derivations $\delta$ such that the property $RIGP_{k(t,\sigma,\delta)/k}$ is satisfied.
\newline 
The following theorem is a restatement in this new context of the results from \cite{BEH}:

\begin{thm} Let $k$ be a field that has a finite-order automorphism $\sigma$ such that the center of $k$ contains an ample field fixed by $\sigma$. Then the property $RIGP_{k(t,\sigma)/k}$ holds. \end{thm}

\textbf{Proof:} Let $\mathbf{k}$ be an ample field fixed by $\sigma$ and contained in the center of $k$, and let $n$ be a multiple of the order of $\sigma$. Let $G$ be a finite group. There exists a commutative field extension $L/\mathbf{k}(t^n)$ that is regular over $\mathbf{k}$ and Galois with Galois group isomorphic to $G$. In \cite{BEH}, it is shown that the $\mathbf{k}(t^n)$-algebra $k(t,\sigma) \otimes_{\mathbf{k}(t^n)} L$ is a field, thus an exterior Galois extension of $k(t,\sigma)$ with Galois group isomorphic to $G$. Furthermore, there exists a $k$-embedding $$k(t,\sigma) \otimes_{\mathbf{k}(t^n)} L \rightarrow k((t^n,\sigma)).$$ Now, $k((t^n,\sigma))/k$ is a regular extension because it is contained in $k((t,\sigma))$, so $k(t,\sigma) \otimes_{\mathbf{k}(t^n)} L$ is regular over $k$, and $RIGP_{k(t,\sigma)/k}$ holds.$\blacksquare$
\newline 
The theorem above is the one we generalize in the subsequent text. The following theorem is a similar statement when the field of fractions is also twisted by a derivation. It is derived from \cite{Des21}, and its proof is essentially the same as that of Theorem 3.

\begin{thm} Let $k$ be a field with center $\mathbf{c}$, $\sigma$ an automorphism of $k$, $\delta$ a $\sigma$-derivation of $k$, and $\mathbf{c}_{cte(\delta)}$ the field of $\delta$-constants of $\mathbf{c}$. If there exists an integer $n \geq 1$ such that $t^n \in Z(k(t,\sigma,\delta))$ and if the subfield $\mathbf{c}_0$ of $\mathbf{c}_{cte(\delta)}$ invariant under $\sigma$ contains an ample field, then the property $\mathrm{RIGP}_{k(t,\sigma,\delta)/k}$ is satisfied. \end{thm}

\textbf{Proof:} Let $\mathbf{k}$ be an ample field contained in $\mathbf{c}_0$ and let $n$ be an integer such that $t^n \in Z(k(t,\sigma,\delta))$. Let $G$ be a finite group. There exists a commutative field extension $L/\mathbf{k}(t^n)$ that is regular over $\mathbf{k}$ and Galois with Galois group isomorphic to $G$. In \cite{Des21}, it is shown that the $\mathbf{k}(t^n)$-algebra $k(t,\sigma,\delta) \otimes_{\mathbf{k}(t^n)} L$ admits a $k$-embedding $k(t,\sigma) \otimes_{\mathbf{k}(t^n)} L \rightarrow k_{\sigma,\delta}((z^n))$. Now, $k_{\sigma,\delta}((z^n))/k$ is a regular extension because it is contained in $k_{\sigma,\delta}((z))$, so $k(t,\sigma,\delta) \otimes_{\mathbf{k}(t^n)} L$ is regular over $k$, and $\mathrm{RIGP}_{k(t,\sigma,\delta)/k}$ holds.
\section{Permanence of the Galois group in towers of Galois extensions}
This section is dedicated to the study of inductive limits of Galois extensions equipped with natural isomorphisms between their Galois groups. It serves as a preliminary section for the study of the regular inverse Galois problem for towers of skew rational fractions fields , which will be addressed in section 4. We start by proving a generalized Artin lemma, then we establish a criterion for the preservation of the Galois group in an inductive limit, and finally, we apply this criterion to inductive limits of skew rational fractions fields.

\subsection{Generalized Artin Lemma}
We present an analogue for arbitrary fields of the version of the Artin lemma for profinite groups stated in \cite{Fr} (lemma 1.3.2).
\newline
Before stating this result, we need to understand the structure of sub-extensions of fields generated by a family of elements. If $R \subset R'$ are rings and $\mathcal{F}$ is a family of elements from $R'$, we denote by $R[\mathcal{F}]$ the subring of $R'$ generated by $R$. This is the set of finite sums of finite products of elements from $R$ and elements from $\mathcal{F}$.
\begin{lem}
Let $h/k$ be a field extension, $I$ a set, and $\mathcal{F}:=(x_i)_{i\in I}$ a family of elements from $h$. Define $R_0 := k[\mathcal{F}]$, and for $n \in \mathbb{N}$, recursively define $R_{n+1} := R_n[R_n^{-1}]$ and $k(\mathcal{F}) := \bigcup_{n \in \mathbb{N}} R_n$.
The following two assertions are true:
\newline
(i) $k(\mathcal{F})$ is the subfield of $h$ generated by $k$ and the family $(x_i)_{i \in I}$.
\newline
(ii) Let $f$ be a ring endomorphism of $h$. If $f(k) \subset k$ and $f(\mathcal{F}) \subset \mathcal{F}$, then the restriction of $f$ to $k(\mathcal{F})$ is an endomorphism of this field. Moreover, this restriction is the identity if and only if $f$ induces the identity on $k$ and on $\mathcal{F}$.
\end{lem}
\textbf{Proof:} By definition, each of the $R_n$ is a subring of $h$. Since $k(\mathcal{F})$ is an increasing union of subrings of $h$, it is also a subring of $h$. Furthermore, for $n \in \mathbb{N}$ and $x \in R_n \setminus \{0\}$, $x^{-1} \in R_{n+1}$, so $k(\mathcal{F})$ is stable under inversion, and hence a subfield of $h$. Since $k(\mathcal{F})$ contains $k$ and the family $\mathcal{F}$, it contains the subfield of $h$ generated by $k$ and $\mathcal{F}$. Additionally, $R_0$ is the subring of $h$ generated by $k$ and the family $\mathcal{F}$, so it is contained in the subfield of $h$ generated by $k$ and the family $(x_i)_{i \in I}$. By induction, every field containing $R_0$ contains $R_n$ for all $n \in \mathbb{N}$, and thus contains $k(\mathcal{F})$. We have thus proven assertion (i).
\newline
For (ii), it is immediate that if $f$ is a ring endomorphism of $h$ satisfying the hypotheses of (ii), then $f(R_0) \subset R_0$. Since $f$ is a ring homomorphism, it follows by induction that for all $n \in \mathbb{N}$, $f(R_n) \subset R_n$, and hence $f(k(\mathcal{F})) \subset k(\mathcal{F})$. Since $k(\mathcal{F})$ is the field of fractions of $h$ over $k[\mathcal{F}]$, the restriction of $f$ to $k(\mathcal{F})$ is the identity.
\newline
\newline
Now that we have this result, we state and prove the announced Artin's Lemma:

\begin{lem}[Generalized Artin's Lemma]
Let $G$ be a profinite group acting faithfully by automorphisms on a field $h$. Suppose that for every $x \in h$, the stabilizer $S(x) := \{\sigma \in G \mid \sigma(x) = x\}$ is an open subgroup of $G$. Let $k := h^G$ be the subfield of $h$ fixed by all elements of $G$. If $h/k$ is an exterior extension, then it is algebraic and Galois, and $\mathrm{Gal}(h/k) = G$. 
\end{lem}
\textbf{Proof:} For $n \in \mathbb{N}$, and $x_1, \ldots, x_n \in h$, the group $H_{x_1, \ldots, x_n} = \bigcap_{i=1}^n S(x_i)$ is open in $G$. Therefore, the intersection $N_{x_1, \ldots, x_n}$ of all conjugates of $H_{x_1, \ldots, x_n}$ is a normal subgroup of finite index in $G$ (see for example \cite{Prof} Lemma 3.1.2). Hence, $G/N$ is a finite group. Let $l := k(Gx_1, \ldots, Gx_n)$ be the subextension of $h/k$ generated by the orbits of the $x_i$ under the action of $G$. By assumption, $G$ acts as the identity on $k$ and $G$ permutes the orbits $Gx_i$ for $1 \leq i \leq n$. Lemma 2 thus allows us to conclude that $G$ acts by automorphisms on $l$.
\newline
The subgroup of $G$ acting as the identity on $l$ is $N_{x_1, \ldots, x_n}$. Indeed, if $\sigma \in G$ induces the identity on $l$, then for all $g \in G$ and $1 \leq i \leq n$, we have $\sigma(g(x_i)) = g(x_i)$, and thus $g^{-1}\sigma g(x_i) = x_i$, meaning $g^{-1}\sigma g \in S(x_i)$. Therefore, $g^{-1}\sigma g \in H_{x_1, \ldots, x_n}$, and $\sigma \in gH_{x_1, \ldots, x_n}g^{-1}$. Thus, $\sigma \in N_{x_1, \ldots, x_n}$. Conversely, if $\sigma \in N_{x_1, \ldots, x_n}$, then for every $g \in G$, we have $g^{-1} \sigma g \in H_{x_1, \ldots, x_n}$, so $\sigma(g(x_i)) = g(x_i)$, and by assumption, $\sigma$ fixes all elements of $k$. Thus, by Lemma 2-(ii), $\sigma$ fixes all elements of $k(x_1, \ldots, x_n) = l$. We have shown that $G/N_{x_1, \ldots, x_n}$ is a finite group acting faithfully by automorphisms on $l$. Moreover, an element $x \in l$ is fixed by $G/N_{x_1, \ldots, x_n}$ if and only if it is fixed by $G$, so the subfield of $l$ fixed by $G/N_{x_1, \ldots, x_n}$ is $k$. Since $h/k$ is an exterior extension, so is $l/k$. Thus, $G/N_{x_1, \ldots, x_n}$, considered as a group of automorphisms of $l$, is an $N$-group of finite cardinality, and the non-commutative version of Artin's lemma for $N$-groups of finite reduced cardinality (\cite{Cohn_1995} Theorem 3.3.7) applies: $l/k$ is a finite Galois exterior extension of degree $[l:k] = |G/N_{x_1, \ldots, x_n}|$. In particular, every element of $h$ is contained in a finite-degree extension of $k$: that is $h/k$ is algebraic.
\newline
$h/k$ is also a Galois extension: indeed, by assumption, $G$ injects into $\mathrm{Aut}(h/k)$ and $h^G = k$, so $k \subset h^{\mathrm{Aut}(h/k)} \subset h^G = K$. Therefore, $h/k$ is an algebraic Galois exterior extension, and the Galois correspondence holds in exactly the same way as for Galois algebraic extensions of commutative fields (\cite{Cohn_1995} Theorem 3.8.4).
\newline
Since the action of $G$ is faithful, the projective system of $G/N_{x_1, \ldots, x_n}$ is cofinal in the projective system of quotients of $G$ by open normal subgroups, and $G \simeq \varprojlim G/N_{x_1, \ldots, x_n}$. On the other hand, $\mathrm{Gal}(h/k) \simeq \varprojlim \mathrm{Gal}(k(Gx_1, \ldots, Gx_n)/k)$, so the isomorphisms between $G/N_{x_1, \ldots, x_n}$ and $\mathrm{Gal}(k(Gx_1, \ldots, Gx_n)/k)$ induce an isomorphism between $G$ and $\mathrm{Gal}(h/k)$. \newline
\newline
\textbf{Remark:} The argument often used for commutative fields, that a field extension $h/k$ which is a union of finite Galois extensions $l_i/k$ is Galois, is not always applicable for skew field extensions: the usual argument for showing that $h{\mathrm{Aut}(h/k)} = k$ uses the fact that if $l_i/k$ is finite and Galois, then for any extension $m/k$, $ml_i/l_i$ is a Galois extension. However, this is proved using the fact that Galois extensions of commutative fields are normal and separable, concepts that rely on relations between finite algebraic field extensions and polynomial rings, which have no known equivalent for skew fields. Thus, it is necessary here to use the action of the group $G$ to conclude that the extension in question is Galois. However, if $h/k$ satisfies an additional lifting property of isomorphisms as considered in (\cite{DesMon}), then it is true that such an extension is Galois.
\subsection{Isomorphisms of Towers of Galois Extensions}

Now that we have an adapted version of Artin's Lemma, we establish that inductive systems of algebraic Galois extensions with the same Galois group $G$, equipped with natural isomorphisms between these Galois groups, give algebraic Galois extensions with group $G$ by passing to the limit.

Let $(I, \leq)$ be a filtered ordered set, and let $(h_i)_{i \in I}$ and $(k_i)_{i \in I}$ be inductive systems of fields. Denote by $\psi_{ij}$ the transition morphisms $h_i \to h_j$, $\phi_{ij}$ the transition morphisms $k_i \to k_j$, and let $h$ and $k$ be the inductive limits of these systems: $h := \varinjlim h_i$, $k := \varinjlim k_i$. Then we have the following result:

\begin{prop} Let $G$ be a profinite group. Suppose that for every $i \in I$, $f_{i} : k_{i} \to h_{i}$ is an algebraic Galois extension, that we have profinite group isomorphisms $\omega_i : G \to \mathrm{Gal}(h_i / k_i)$ for every $i \in I$, and that the following diagrams commute for $i \leq j \in I$:

\begin{equation}
\begin{tikzcd}
k_{j} \arrow[rr, "f_{j}"]                  &  & h_{j}                  \\
                                    &  &                          \\
k_i \arrow[rr, "f_i"] \arrow[uu, "\phi_{ij}"] &  & h_i \arrow[uu, "\psi_{ij}"]
\end{tikzcd}
\end{equation}

Moreover, assume that for all $i \leq j \in I$, every $\sigma \in G$ and $x \in h_i$:

\begin{equation} \psi_{ij}(\omega_i(\sigma)(x)) = \omega_j(\sigma)(\psi_{ij}(x)). \end{equation}

Then $h/k$ is an algebraic Galois extension, and there exists an isomorphism $\omega : G \to \mathrm{Gal}(h/k)$ such that for every $i \in I$, $x \in h_i \subset h$, and $\sigma \in G$, we have $\omega(\sigma)(x) = \omega_i(\sigma)(x)$. \end{prop}

\textbf{Proof:} Define a morphism of fields $f : k \to h$. For every $i \in I$, the composition of the morphism $f_i : k_i \to h_i$ and the canonical inclusion morphism $h_i \to h$ gives a morphism which we denote by $f_i : k_i \to h$. Since diagram (1) commutes for all $i \leq j \in I$, the $f_i : k_i \to h$ form a compatible system with the inductive system $(k_i)_{i \in I}$, and thus define a morphism $f : k \to h$ such that for every $i \in I$ and $x \in k_i$, we have $f(x) = f_i(x)$. Hence, $f$ defines a morphism (non-zero, since the $f_i : k_i \to h$ are non-zero) $k \to h$, that is, a field extension $h/k$. From now on, we denote by $x \in k \subset h$ the element $f(x)$ in $h$.

Let $g$ be an inner automorphism of $h$ such that $f|_k = \mathrm{Id}_k$. Then there exists an element $\alpha \in h$ such that $g = I(\alpha)$. Let $i \in I$ be such that $\alpha \in h_i \subset h$. Then for all $j \geq i$, $\alpha \in h_j \subset h$ and $I(\alpha)$ is an outer automorphism of $h_j$. But since $g$ is the identity on $k$, $g$ restricted to $k_j \subset h_j$ is also the identity of $k_j$. Hence, $g|_{h_j} \in \mathrm{Gal}(h_j / k_j)$. Since it is an inner automorphism of $h_j$ and $h_j / k_j$ is an outer extension, the restriction of $g$ to $h_j$ is the identity. Thus, $g$ induces the identity on each $h_i$ for $i \in I$, and consequently $g = \mathrm{Id}_h$. Therefore, $h/k$ is an outer extension.

We claim that for $\sigma \in G$, the following formula defines a $k$-automorphism of $h$:

\begin{align*} g_{\sigma} : h & \to h \\ h_i \ni x & \mapsto \omega_i(\sigma)(x). \end{align*}

Condition (2) ensures that $\omega_i(\sigma)(x)$ depends only on the class of $x$ in $h$, so $g_\sigma$ is a well-defined function. If $x, y \in h_i$ for some $i \in \mathbb{N}$, we have
 $$g_{\sigma}(x-y)=\omega_i(\sigma)(x-y)=\omega_i(\sigma)(x)-\omega_i(\sigma)(y)=g_{\sigma}(x)-g_{\sigma}(y)$$ and $$g_{\sigma} (xy^{-1})=\omega_i(\sigma)(xy^{-1})=\omega_i(\sigma)(x)\omega_i(\sigma)(y)^{-1}=g_{\sigma}(x)g_{\sigma}(y)^{-1}.$$ 
Since any pair of elements of $h$ is contained in some $h_i \subset h$, it follows that $g_\sigma$ is a field morphism. Let $x \in h$. For every $i \in \mathbb{N}$ such that $x \in h_i \subset h$, since $\omega_i(\sigma)$ is an automorphism of $h_i$, there exists a unique $y \in h_i$ such that $\omega_i(\sigma)(y) = x$. Condition (2) implies that the class of $y$ in $h$ does not depend on the chosen integer $i$, so there exists a unique $y \in h$ such that $g_{\sigma}(y) = x$. Thus, $g_{\sigma}$ is an automorphism of $h$. If $\sigma \neq \sigma'$ and $i$ is an integer, then $g_\sigma|_{h_i} = \omega_i(\sigma) \neq \omega_i(\sigma') = g_{\sigma'}|_{h_i}$, so $g_{\sigma} \neq g_{\sigma'}$. Since for every $i \in I$, the map $\omega_i$ is a group morphism, for all $x \in h$, and $\sigma, \sigma' \in G$, we have
$$g_{\sigma^{-1} \sigma'}(x)= g_{\sigma}^{-1}(x) g_{\sigma'}(x).$$

Moreover, for $i \in I$, $\sigma \in G$, and $x \in k_i$, we have $g_\sigma(x) = \omega_i(\sigma)(x) = x$. Hence, $g_\sigma$ is a $k$-automorphism of $h$. Therefore, $\omega : G \to \mathrm{Aut}_k(h), \sigma \mapsto g_\sigma$ is an injective group morphism.

Now, we show that $h^G = k$. We have already established that $k \subset h^G$. Let $x \in h \setminus k$. Let $i \in I$ such that $x \in h_i \subset h$. Then, since $h_i^G = k_i$, there exists $\sigma \in G$ such that $\omega_i(\sigma)(x) \neq x$. In particular, $\omega(\sigma)(x) \neq x$, so $x \notin h^G$. Therefore, $h^G \subset k$, and we conclude that $h^G = k$.

Finally, we apply the generalized Artin's Lemma. $G$ is a profinite group acting faithfully by automorphisms on $h$. Let $x \in h$, and let $i \in I$ such that $x \in h_i \subset h$. Then $k_i(x)/k_i$ is a finite subextension of $h_i/k_i$, and $S(x)$ is the subgroup of $\mathrm{Gal}(h_i/k_i)$ that induces the identity on $k_i(x)$. The Galois correspondence for outer algebraic extensions ensures that $S(x)$ is an open subgroup of $G$. Moreover, since $h^G = k$, it follows that $h/h^G$ is an outer extension. Thus, the hypotheses of the generalized Artin's Lemma are satisfied, and we conclude that $h/k$ is an outer Galois extension with Galois group isomorphic to $G$.

\textbf{Remarks:} a) If the $h_i$ and $k_i$ were commutative fields, this proposition could be proven without using Artin's Lemma. In the non-commutative case, it remains to verify that the extension is outer and that the considered subgroup of the automorphism group of $h$ isomorphic to $G$, is an $N$-group of automorphisms of $h$. More precisely, if the $h_i$ and $k_i$ are commutative, saying that the extension $h/k$ is an algebraic Galois extension is equivalent to saying that it is algebraic, normal, and separable. However, the notions of normality and separability depend on the properties of polynomial evaluation morphisms, for which no suitable equivalent is known for skew fields.

b) If the group $G$ is finite, the hypotheses of Artin's non-commutative Lemma for groups of automorphisms of a field with finite reduced cardinality are satisfied, so the hypothesis that the $h_i/k_i$ are algebraic extensions is always satisfied. However, if $G$ is profinite but not finite, then the hypothesis that the stabilizer of every element in $h$ is an open subgroup of $G$ is needed. If $G$ acts faithfully and the extensions $h_i/k_i$ are Galois extensions with group $G$, with the isomorphisms between Galois groups and $G$ satisfying the hypotheses of the proposition, it is equivalent to say that for every $x$, the stabilizer of $x$ in $G$ is open and that $x$ is algebraic over each of the $k_i$ and over $k$.

\textbf{Examples:} a) A first example is the commutative case. Let $k_n$ be an increasing sequence of number fields such that $-1$ is not a square in any $k_n$. The proposition then says that if $k = \varinjlim k_n$, $k(\sqrt{-1}) / k$ is a Galois extension with Galois group $\mathbb{Z}/2\mathbb{Z}$, where the non-trivial element is $a + b\sqrt{-1} \mapsto a - b\sqrt{-1}$.

b) For a non-commutative example, consider $m_n$ an increasing sequence of subfields of $\mathbb{C}$, algebraic over $\mathbb{Q} \subset \mathbb{C}$, not contained in $\mathbb{R} \subset \mathbb{C}$, and stable under conjugation. Let $f\subset \mathbb{R}$ be a finite Galois extension of $\mathbb{Q}$ that is linearly disjoint from the $m_n$, so that $fm_n/m_n$ is a finite Galois extension with Galois group $\mathrm{Gal}(f/\mathbb{Q})$ for every natural number $n$. Let $\tau : \mathbb{C} \to \mathbb{C}$ be complex conjugation. We can then set $k_i := m_i(t, \tau)$, $h_i := fm_i(t, \tau)$, and all the field morphisms used in the proposition will be inclusions. Then the proposition shows that $\varinjlim fm_n(t, \tau) / \varinjlim m_n(t, \tau)$ is a Galois extension with Galois group $\mathrm{Gal}(f/\mathbb{Q})$.

The extension $h_i/k_i$ is an outer Galois extension with Galois group $\mathrm{Gal}(f/\mathbb{Q})$. Indeed, if we take $(f_\sigma)_{\sigma \in \mathrm{Gal}(f/\mathbb{Q})}$ as a normal basis of $f/\mathbb{Q}$, it is also a normal basis of $fm_n/m_n$. Therefore, if $x \in fm_n$, $x = \sum a\sigma f_\sigma$ with $\sigma \in \mathrm{Gal}(f/\mathbb{Q})$ and $a_\sigma \in m_n$, then $\tau\sigma(x) = \tau(\sum a_{\gamma} f_{\sigma\gamma}) = \sum \tau(a_{\gamma}) f_{\sigma\gamma}$ and $\sigma\tau(x) = \sigma(\sum \tau(a_{\gamma}) f_{\gamma}) = \sum \tau(a_{\gamma}) f_{\sigma\gamma}$. This commutation of $\sigma$ and $\tau$ allows us to define an automorphism $\sigma$: $fm_n[t,\tau] \rightarrow fm_n[t,\tau], \sum a_n t^n \mapsto \sum \sigma(a_n) t^n$. This in turn defines an automorphism of $fm_n(t,\tau)$. We have thus defined a subgroup of the automorphism group of the field $fm_n(t,\tau)$ that fixes $m_n(t,\tau)$. Each of these automorphisms (except for the identity) is outer since they act non-trivially on the center (because the center of $\mathbb{C}[t,\tau]$ is $\mathbb{R}[t,\tau]$ and $f[t,\tau] \subset \mathbb{R}[t,\tau]$). This group $G$ of automorphisms is therefore an $N$-group, which is canonically identified with $\mathrm{Gal}(f/\mathbb{Q})$, and $fm_n(t,\tau)/fm_n(t,\tau)^G$ is an outer Galois extension with Galois group $\mathrm{Gal}(f/\mathbb{Q})$. We also have $[fm_n(t,\tau): fm_n(t,\tau)^G] = |G|$. Now, $[fm_n(t,\tau): m_n(t,\tau)] = [fm_n: m_n] = |G|$ (for the proof of the first equality, see \cite{Des24}, Theorem 2), so $fm_n(t,\tau)^G = m_n(t,\tau)$ and $fm_n(t,\tau)/m_n(t,\tau)$ is an outer Galois extension with Galois group $\mathrm{Gal}(f/\mathbb{Q})$. Therefore, the hypotheses of the proposition are verified.
\subsection{Skew Fraction Fields Extensions}
We apply the results of the previous section when the fields $(k_i)_{i \in I}$ in diagram (1) are twisted fraction fields, and the fields $(h_i)_{i \in I}$ in the same diagram are obtained from the $k_i$ by scalar extensions.
\newline
Let $h/k$ be a field extension such that there exists an increasing sequence $(h_i)_{i \in \mathbb{N}}$ of subfields of $h$ with $h_i/k$ being an outer Galois extension of finite degree $m_i$. Let $\sigma$ be an automorphism of $h$ such that the restriction of $\sigma$ to $k$ is the identity, and for every natural number $i$, the restriction of $\sigma$ to $h_i$ is an automorphism of $h_i$, which we will still denote as $\sigma$ when there is no ambiguity. We are interested in the Galois extensions of the field $h(t,\sigma)$. Suppose $\sigma = \mathrm{Id}_h$, $\mathbf{k} \subset Z(k)$, and there exists $\mathbf{k}(t) \subset l \subset \mathbf{k}((t))$ such that $L/\mathbf{k}(t)$ is a finite Galois extension with Galois group $G$. Then the ring $h_i(t,\sigma) \otimes_{\mathbf{k}(t)} L$ is a field, and $h_i(t,\sigma) \otimes_{\mathbf{k}(t)} L / h_i(t,\sigma)$ is a Galois extension with Galois group $G$. Each of the following commutative squares satisfies the assumptions of Proposition 4:

$$\begin{tikzcd} h_{i+1}(t,\sigma) \arrow[rr] & & h_{i+1}(t,\sigma) \otimes_{\mathbf{k}(t)} L \\ & & \\ h_i(t,\sigma) \arrow[rr] \arrow[uu] & & h_i(t,\sigma) \otimes_{\mathbf{k}(t)} L \arrow[uu] 
\end{tikzcd}$$

As a result, $h(t,\sigma) \otimes_{\mathbf{k}(t)} L = \varinjlim h_i(t,\sigma) \otimes_{\mathbf{k}(t)} L / \varinjlim h_i(t,\sigma) = h(t,\sigma)$ is an outer Galois extension with Galois group $G$.
\newline
We wish to apply this construction to the case where the indeterminate $t$ is not central, i.e., when $\sigma$ is not the identity. Fix $i \in \mathbb{N}$ and assume $\sigma|_{h_i} \neq \mathrm{Id}_{h_i}$. The multiplication by elements of $\mathbf{k}(t) \subset h_i(t,\sigma)$ does not give $h_i(t,\sigma)$ a structure of $\mathbf{k}(t)$-algebra. Indeed, if $a$ and $b$ are elements of $h_i$ and $a \neq \sigma(a)$, we have: 
 $$(ta)(tb)=\sigma (a)\sigma^2 (b)t^2\neq \sigma^2 (a)\sigma^2(b) t^2=t^2 (ab).$$
As a result, the $\mathbf{k}(t)$-bilinear map on $h_i(t,\sigma) \otimes_{\mathbf{k}(t)} L$: $$(a \otimes b, c \otimes d) \mapsto ab \otimes cd$$ does not, in general, define a ring structure on $h_i(t,\sigma) \otimes_{\mathbf{k}(t)} L$.  However, since $L/\mathbf{k}(t)$ is separable, the primitive element theorem ensures the existence of an element $a \in L$ such that $L = \mathbf{k}(t)(a)$. Since $L \subset \mathbf{k}((t))$, this element $a$ can be written as $a = S(t) = \sum_{n \in \mathbb{Z}} S_n t^n$. We then denote $S(t^m)$ as the Laurent series $\sum_{n \in \mathbb{Z}} S_n t^{mn}$ and define $L_m := \mathbf{k}(t^m)(S(t^m))$. We have $\mathbf{k}(t^m) \subset L_m \subset \mathbf{k}((t^m))$ and $L_m/\mathbf{k}(t^m)$ is a Galois extension with Galois group $G$ (see \cite{BEH} Proposition 2.1.2). Since $\sigma|_{h_i}$ is in the Galois group of $h_i/k$, the order of $\sigma$ divides $m_i$. Therefore, $t^{m_i}$ is central in $h_i(t,\sigma)$, and the field $h_i(t,\sigma)$ is a $\mathbf{k}(t^{m_i})$-algebra. The $\mathbf{k}(t^{m_i})$-algebra $h_i(t,\sigma) \otimes_{\mathbf{k}(t^{m_i})} L_{m_i}$ is an outer Galois extension of $h_i(t,\sigma)$ with Galois group $G$. Our goal is to adapt the construction for the case where $t$ is central by filtering $h$ through sub-extensions on which $\sigma$ has finite order, so that we can use Proposition 4. To do this, we first prove a lemma that establishes sufficient conditions to obtain a commutative diagram as in Proposition 4 from two extensions $k(t,\sigma^b) \rightarrow k(t,\sigma^b) \otimes_{\mathbf{k}(t^n)} L_n$ and $h(t,\sigma^a) \rightarrow h(t,\sigma^a) \otimes_{\mathbf{k}(t^{nm})} L_{nm}$.
\newline
\begin{lem}
Let $h/k$ be a field extension, $\mathbf{k}$ a subfield of $Z(h)$, the center of $h$, and $\sigma$ an automorphism of $h$ that induces an automorphism of $k$ by restriction. Let $a, b, n, m$ be natural numbers. Suppose that $ord(\sigma|_k)$ divides $bn$, $ord(\sigma)$ divides $anm$, and $am \equiv b \pmod{ord(\sigma|_k)}$. We define a map:
\begin{align*}
f : k[t,\sigma^b]&\rightarrow h[t,\sigma^a] \\
\sum_{i=0}^l a_i t^i&\mapsto\sum_{i=0}^l a_i t^{mi}
\end{align*}
Then $f$ extends to a morphism of $\mathbf{k}(t^n)$-algebras $k(t, \sigma^b) \rightarrow h(t, \sigma^a)$, and if we denote by $g_{n, nm} : L_n \rightarrow L_{nm}$ the isomorphism defined by $g_{n, nm}(S(t^n)) = S(t^{nm})$, the following diagram is a commutative diagram of $\mathbf{k}(t^n)$-algebras:

$$\begin{tikzcd} {h(t, \sigma^a)} \arrow[rrrr, "x \mapsto x \otimes 1"] & & & & {h(t, \sigma^a) \otimes_{\mathbf{k}(t^{nm})} L_{nm}} \\ & & & & \\ {k(t, \sigma^b)} \arrow[rrrr, "x \mapsto x \otimes 1"'] \arrow[uu, "f"'] & & & & {k(t, \sigma^b) \otimes_{\mathbf{k}(t^n)} L_n} \arrow[uu, "f \otimes g_{n, nm}"'] \end{tikzcd}$$
\end{lem} \textbf{Proof:} The horizontal arrows are ring morphisms as the multiplication in $k(t, \sigma^b)$ induces a structure of $\mathbf{k}(t^n)$-algebra and the multiplication in $h(t, \sigma^a)$ induces a structure of $\mathbf{k}(t^{nm})$-algebra. This holds if and only if $t^n$ is in the center of $k(t, \sigma^b)$ and $t^{nm}$ is in the center of $h(t, \sigma^a)$. Since for every $\lambda \in k$, $t^n \lambda = \sigma^{bn}(\lambda) t^n$ in $k(t, \sigma)$, and for every $\lambda \in h$, $t^{nm} \lambda = \sigma^{anm}(\lambda) t^{nm}$ in $h(t, \sigma^a)$, $t^n$ is in the center of $k(t, \sigma^b)$ if and only if $\mathrm{ord}(\sigma|k)$ divides $bn$, and $t^{nm}$ is in the center of $h(t, \sigma^a)$ if and only if $\mathrm{ord}(\sigma)$ divides $anm$. \newline
The way $f$ is defined on twisted polynomial rings ensures that for every $\lambda \in k$ and $k \in \mathbb{N}$, $f(\lambda t^k) = \lambda t^{km} = f(\lambda) f(t^k)$. Moreover, $f(1) = 1$ and for all $x, y \in k[t, \sigma^b]$, we have $f(x + y) = f(x) + f(y)$. We now need to verify that $f(t \lambda) = f(t) f(\lambda)$. But $$f(t \lambda) = f(\sigma^b(\lambda) t) = \sigma^b(\lambda) t^m$$ and $$f(t) f(\lambda) = t^m \lambda = \sigma^{am}(\lambda) t^m.$$ Therefore, we obtain a ring homomorphism if and only if for all $\lambda \in k$, we have $\sigma^b(\lambda) = \sigma^{am}(\lambda)$, i.e., if $am \equiv b \pmod{\mathrm{ord}(\sigma|k)}$. This is an injective ring homomorphism between Ore domains, and by the properties of Ore domains, it induces a field homomorphism between their fraction fields, $f : k(t, \sigma^b) \rightarrow h(t, \sigma^a)$.
\newline
The right arrow is defined as follows: $h(t, \sigma^a)$ and $L_{nm}$ have a structure of $\mathbf{k}(t^n)$-vector spaces induced by the field isomorphism \begin{align*}\mathbf{k}(t^n) &\rightarrow \mathbf{k}(t^{nm})\\ P(t^n)& \mapsto P(t^{nm}),\end{align*} and the map \begin{align*}k(t, \sigma^b) \times L_n &\rightarrow h(t, \sigma^a) \otimes_{\mathbf{k}(t^{nm})} L_{nm}\\ (x, \sum_{i=0}^r Q_i(t^n) S(t^n)^i) &\mapsto f(x) \otimes (\sum_{i=0}^r Q_i(t^{nm}) S(t^{nm})^i)\end{align*} is $\mathbf{k}(t^n)$-bilinear for this vector space structure. Therefore, $x \otimes S \mapsto f(x) \otimes g_{n, nm}(S)$ defines a morphism of $\mathbf{k}(t^n)$-vector spaces. In particular, it is a morphism of abelian groups. Since $f : k(t, \sigma^b) \rightarrow h(t, \sigma^a)$ and $g_{n, nm} : L_n \rightarrow L_{nm}$ are also ring homomorphisms, it is even a morphism of $\mathbf{k}(t^n)$-algebras, and in particular, a ring homomorphism.
\newline
\textbf{Remarks:} (i) Under the assumptions of the lemma, or even with a finite number of extensions, one could simply extend scalars by taking a tensor product of $\mathbf{k}(t^N)$-algebras with some $N$ such that the multiplication in each of the fields induces a structure of $\mathbf{k}(t^N)$-algebra. However, if one wants to use the lemma to apply Proposition 1 to an inductive limit of $h_i(t, \sigma)$, there will be infinitely many such diagrams, and if the order of $\sigma$ in the inductive limit is infinite, there will be no integers $N$ such that the $h_i(t, \sigma)$ are all $\mathbf{k}(t^N)$-algebras.
\newline
(ii) It is important to emphasize the necessity that $t^n$ (respectively $t^{nm}$) be central in $k(t, \sigma^b)$ (respectively $h(t, \sigma^a)$). If these elements are not central, $k(t, \sigma^b)$ and $h(t, \sigma^a)$ will no longer be $\mathbf{k}(t^n)$ or $\mathbf{k}(t^{nm})$-algebras, and one would need to check whether the formula used to define $g$ does indeed define a homomorphism of (non-commutative) rings. One might also wonder whether, given two field homomorphisms $k \rightarrow h$ and $k \rightarrow L$, one could define a $k$-ring $h \otimes_k L$ to work with, but then one faces the fact that the tensor product in the category of $(k,k)$-bimodules only takes into account the left module structure on $h$ and the right module structure on $L$. The ring structure on this tensor product defined by $(a \otimes b)(c \otimes d) = (ab \otimes cd)$ would make the elements of $h$ and $L$ commute: indeed, we would have $(1 \otimes b)(a \otimes 1) = (a \otimes 1)(1 \otimes b) = a \otimes b$.
\section{Application to the Inverse Galois Problem} We apply the results from the previous section to show that when $h$ satisfies sufficient arithmetic properties, the property $\mathrm{RIGP}_{H/h}$ is verified when $H$ is an inductive limit of twisted fraction fields, constructed using the mentioned arithmetic properties and whose union of coefficient fields is $h$.  

\subsection{Tower with a Fixed Galois Group} We extend the results of \cite{BEH} to towers of skew fraction fields using the results from the previous section. \newline The following theorem provides arithmetic conditions for constructing extensions of skew fraction fields towers with a given Galois group $G$ from a pair $(h, \sigma)$, where $h$ is a field and $\sigma \in \mathrm{Aut}(h)$ is locally of finite order. 

\begin{thm}
Let $h/k$ be a field extension such that there exists a sequence $h_n$ of subfields of $h$ with $h_n \subset h_{n+1}$ for all $n \in \mathbb{N}$, $h = \varinjlim h_n$, and for every natural number $n$, $h_n/k$ is an outer Galois extension of degree $d_n$ with $Gal(h_n/k)$ cyclic and generated by an element $\sigma_n$ such that for $m \leq n$, $\sigma_n|_{h_m} = \sigma_m$. Let $\mathbf{k} \subset Z(k)$ be a commutative field, and let $L/\mathbf{k}(t)$ be a Galois extension such that $\mathbf{k}(t) \subset L \subset \mathbf{k}((t))$ and $\mathrm{Gal}(L/\mathbf{k}(t)) \cong G$. \newline Let $(a_n)_{n\in\mathbb{N}}, (l_n)_{n\in\mathbb{N}}$ be sequences of natural numbers such that: 

\begin{align*}
&l_n | l_{n+1} \\
&d_n | \frac{a_{n+1}l_{n+1}}{l_n} - a_n \\
&d_n | l_n a_n
\end{align*}

We have families of morphisms for $n \in \mathbb{N}$ given by 
$$u_n : h_n(t, \sigma^{a_n}) \rightarrow h_{n+1}(t, \sigma^{a_{n+1}}), P(t) \mapsto P(t^{\frac{l_{n+1}}{l_n}})$$
and 
$$v_n : L_{l_n} \rightarrow L_{l_{n+1}}, P(t^{l_n}) \mapsto P(t^{l_{n+1}}).$$

Let $H := \varinjlim h_n(t, \sigma^{a_n})$ and $H_G := \varinjlim h_{n}(t, \sigma_n^{a_n}) \otimes_{\mathbf{k}(t^{k_n})} L_{l_n}$, where the inductive limits are taken according to the transition morphisms induced by $u_n$ and $v_n$. Then there exists a field morphism $H \rightarrow H_G$ such that $H_G/H$ is an outer Galois extension with Galois group isomorphic to $G$.

\end{thm}

\textbf{Proof:} The divisibility assumptions of the theorem ensure that the following diagram satisfies the conditions of Lemma 4 for all $n \in \mathbb{N}$.

\[
\begin{tikzcd}
{h_{n+1}(t, \sigma^{a_{n+1}})} \arrow[rrrr, "x \mapsto x \otimes 1"] & & & & {h_{n+1}(t, \sigma^{a_{n+1}}) \otimes_{\mathbf{k}(t^{l_{n+1}})} L_{l_{n+1}}} \\
& & & & \\
{h_n(t, \sigma^{a_n})} \arrow[rrrr, "x \mapsto x \otimes 1"] \arrow[uu, "u_n"'] & & & & {h_n(t, \sigma^{a_n}) \otimes_{\mathbf{k}(t^{l_n})} L_{l_n}} \arrow[uu, "u_n \otimes v_n"]
\end{tikzcd}
\]

Thus, by Lemma 4, it follows that these arrows are ring homomorphisms. By Proposition 2.1.2 of \cite{BEH}, $h_n(t, \sigma^{a_n}) \otimes_{\mathbf{k}(t^{l_n})} L_{l_n}$ is a field for all $n \in \mathbb{N}$, so the arrows in the previous diagrams are field homomorphisms (non-zero). By Lemma 2.1.1 of \cite{BEH}, the horizontal arrows are outer Galois extensions, and there exists a group isomorphism given by:

\begin{align*}
\Psi : \mathrm{Gal}(L_{l_n}/\mathbf{k}(t^{l_n})) & \rightarrow \mathrm{Gal}(h_n(t, \sigma^{a_n}) \otimes_{\mathbf{k}(t^{l_n})} L_{l_n} / h_n(t, \sigma^{a_n})) \\
\varphi & \mapsto \Psi (\varphi) : h_n(t, \sigma^{a_n}) \otimes_{\mathbf{k}(t^{l_n})} L_{l_n} \rightarrow h_n(t, \sigma^{a_n}) \otimes_{\mathbf{k}(t^{l_n})} L_{l_n} \\
& y \otimes x \mapsto y \otimes \varphi(x).
\end{align*}

Let $\varphi$ be an element of $\mathrm{Gal}(L/\mathbf{k}(t))$, and let $x$ be an element of $L$. There exists a unique sequence $(\varphi_i(x))_{i \in \mathbb{Z}}$ of elements of $\mathbf{k}$ such that $\varphi(x) = \sum_{i=-\infty}^{+\infty} \varphi_i(x) t^i.$ Let $n$ be a non-zero natural number, $y := \sum_{i=-\infty}^{+\infty} y_i t^{ni}$. We define:

$$g_n (\varphi)(y) := \sum_{i=-\infty}^{+\infty} \varphi_i\left( \sum_{i=-\infty}^{+\infty} y_i t^i \right) t^{ni}.$$

Then $g_n$ is a group isomorphism $\mathrm{Gal}(L/\mathbf{k}(t)) \rightarrow \mathrm{Gal}(L_n/\mathbf{k}(t^n))$. Let us take:

$$\omega_n := \Psi \circ g_{l_n}, \mathrm{Gal}(L/\mathbf{k}(t)) \rightarrow \mathrm{Gal}(h_n(t, \sigma^{a_n}) \otimes_{\mathbf{k}(t^{l_n})} L_{l_n} / h_n(t, \sigma^{a_n})).$$

By construction, the family $(\omega_n)_{n \in \mathbb{N}}$ is a family of isomorphisms that satisfies the hypotheses of Proposition 4. Thus, Proposition 4 ensures that the previous diagram defines an embedding $H \rightarrow H_G$ that makes $H_G$ an outer Galois extension of $H$ with Galois group $\mathrm{Gal}(L/\mathbf{k}(t)) \cong G$.
\newline
\textbf{Examples:} (a) Let $p$ be a prime number, $(p_n)_{n \in \mathbb{N}}$ an enumeration of the prime numbers, and let $\Pi_n := \prod_{i=0}^n p_n$. Take $h_n := \mathbb{F}_{p^{\Pi_n}}$ and for every $n \in \mathbb{N}$, let $\sigma$ be the Frobenius automorphism $h_n \rightarrow h_n, x \mapsto x^p$. Let $\mathbb{F}_p(t) \subset L \subset \mathbb{F}_p((t))$ such that $L/\mathbb{F}_p(t)$ is a Galois extension. Taking $m_n$ to be a sequence of integers such that $m_n p_{n+1} \equiv 1 \ (\mathrm{mod} \ \Pi_n)$, $a_n = 1$, $l_n = \alpha_n \Pi_n$ with $\alpha_0 = 1$ and $\alpha_{n+1} = m_n \alpha_n$, the theorem applies: 

$$ \varinjlim \mathbb{F}_{p^{\Pi_n}}(t, \sigma) \otimes_{\mathbb{F}_p(t^{l_n})} L_{l_n} / \varinjlim \mathbb{F}_{p^{\Pi_n}}(t, \sigma) $$

is a Galois extension with Galois group $\mathrm{Gal}(L/\mathbb{F}_p(t))$. 
\newline \newline
(b) Similarly, if $G$ is a group that is the Galois group of a regular extension $L/\mathbb{Q}(X)$, let $k_n$ be a cyclic extension of degree $p_n$ of $\mathbb{Q}$, take $h_n$ as the compositum of the $k_i$ for $i \leq n$, and let $h$ be the compositum of all the $h_n$ for $n \in \mathbb{N}$. If $\sigma$ is a topological generator of $\mathrm{Gal}(h/\mathbb{Q})$, that is, an automorphism acting on each $k_n$ as a generator of the Galois group $\mathrm{Gal}(k_n/\mathbb{Q})$, and we choose $l_n$ as in the previous example, then the theorem applies: 

$$ \varinjlim h_n(t, \sigma) \otimes_{\mathbb{Q}(t^{l_n})} L_{l_n} / \varinjlim h_n(t, \sigma) $$

is a Galois field extension with Galois group $G$.
\subsection{Towers Satisfying $\mathrm{IGP}$}

In this subsection, we apply Theorem 5 to construct twisted fraction field towers that satisfy $\mathrm{IGP}$.

The previous examples give the simplest possible application of the theorem: when we can choose $a_n=1$ for all $n$. The divisibility conditions then simplify as follows: $l_n | l_{n+1}$, $d_n | \frac{l_{n+1}}{l_n} - 1$, and $d_n | l_n$. Suppose that the last condition is satisfied and write $l_n = \alpha_n d_n$. We now need to verify that we have $\alpha_n d_n | \alpha_n d_{n+1}$ and $d_n | \frac{\alpha_{n+1}}{\alpha_n} \frac{d_{n+1}}{d_n} - 1$. The first condition is easily verified if $\alpha_n | \alpha_{n+1}$. In that case, verifying the last condition amounts to finding an integer $m$ such that $m \frac{d_{n+1}}{d_n} \equiv 1 \ (\mathrm{mod} \ d_n)$. Such an integer exists if and only if $d_n$ and $\frac{d_{n+1}}{d_n}$ are coprime. Therefore, if this condition on the degrees $d_n$ is satisfied, we obtain the following theorem:

\begin{thm} Let $h/k$ be a field extension such that there exists a sequence $h_n$ of subfields of $h$ with $h_n \subset h_{n+1}$ for all $n \in \mathbb{N}$, $h = \varinjlim h_n$, and for each natural number $n$, $h_n/k$ is an exterior Galois extension of degree $d_n$ with $Gal(h_n/k)$ cyclic, generated by an element $\sigma_n$ such that for $m \leq n$, $\sigma_n |_{h_m} = \sigma_m$. Suppose that for all natural numbers $n$, $\frac{d_{n+1}}{d_n}$ and $d_n$ are coprime. Then, using the notation from Theorem 3, there exists a sequence $(l_n)_{n \in \mathbb{N}}$ such that if $(a_n)_{n \in \mathbb{N}}$ is the constant sequence $1$, the sequences $(l_n)_{n \in \mathbb{N}}$ and $(a_n)_{n \in \mathbb{N}}$ satisfy the assumptions of Theorem 3, and for any such sequence, $\mathrm{SIGP}_\mathbf{k} \Rightarrow \mathrm{IGP}_{h}$. \end{thm}

\textbf{Remark:} There is no reason for $\varinjlim h_n(t,\sigma)$ to be isomorphic to $h(t,\sigma)$ because the morphisms $h_n(t,\sigma) \rightarrow h_{n+1}(t,\sigma)$ are not the canonical inclusions. Since in this situation $l_{n+1}/l_n$ is never equal to 1, the element $t \in h_n(t,\sigma) \subset H$ has $m$-th roots for arbitrarily large $m$, and thus in particular, there cannot exist a field isomorphism between $H$ and $h(t,\sigma)$ that sends one of these elements to $t \in h(t,\sigma)$. Moreover, Theorem 5 can therefore a priori be used to construct a Galois group over many different fields. Indeed, a single extension $h/k$ may admit several different filtrations $h = \varinjlim h_n$, inducing a different sequence $d_n$: such a choice is equivalent to choosing a decreasing sequence of distinguished open subgroups containing $\sigma$, which forms a neighborhood basis of $1$ in $\mathrm{Gal}(h/k)$. Given a sequence $d_n$ provided by such a filtration, there may exist several sequences $a_n$ and $l_n$ satisfying the divisibility conditions of the theorem. The data of a filtration $(h_n)_{n \in \mathbb{N}}$, along with sequences $(a_n)_{n \in \mathbb{N}}$ and $(l_n)_{n \in \mathbb{N}}$, fixes $H(h_n,a_n,l_n) := \varinjlim h_n(t,\sigma)$, but distinct data give a priori non-isomorphic extensions of $k(t)$ or $h$.
\begin{cor} Let $h/k$ be an exterior Galois extension with $\mathrm{Gal}(h/k)$ procyclic, and such that the supernatural number $p^\infty$ does not divide the supernatural order $\#\mathrm{Gal}(h/k)$ for any prime number $p$. Let $(p_n)_{n \in \mathbb{N}^*}$ be an enumeration of prime numbers, and $\#\mathrm{Gal}(h/k) = \prod_{n \in \mathbb{N}^*} p_n^{l_n}$. Then, for any commutative subfield $\mathbf{k}$ of $k$, $\mathrm{SIGP}_\mathbf{k} \Rightarrow \mathrm{IGP}_H$ with $H := \varinjlim h_n(t,\sigma)$ for $\sigma$ a topological generator of $\mathrm{Gal}(h/k)$, $h_n$ the subfield fixed by $\sigma^{\prod_{1 \leq i \leq n} p^{l_n}}$, and the morphism $h_n(t,\sigma) \rightarrow h_{n+1}(t,\sigma)$ given by $t \mapsto t^{m_n p_{n+1}^{l_{n+1}}}$ for $m_n$ an integer such that $m_n p_{n+1}^{l_{n+1}} \equiv 1 \ (\mathrm{mod} \ \prod_{1 \leq i \leq n} p_i^{l_i})$. \end{cor}

This corollary particularly applies in the case where the center of $k$ contains an ample field:

\begin{cor} Let $\mathbf{k}$ be a commutative ample field, $k$ a field such that $\mathbf{k} \subset Z(k)$, $h/k$ an exterior Galois extension such that $\mathrm{Gal}(h/k)$ is procyclic, and such that the supernatural number $p^\infty$ does not divide the supernatural order $\#\mathrm{Gal}(h/k)$ for any prime number $p$. Then, the property $\mathrm{IGP}_H$ holds for $H$, $\sigma$, and the morphism $h_n(t,\sigma) \rightarrow h_{n+1}(t,\sigma)$ as described in Corollary 1. \end{cor}

\textbf{Remark:} The assumptions can be easily verified if the Galois group is a product of cyclic groups: if $h/k$ is an exterior Galois extension with $\mathrm{Gal}(h/k) = \prod_{i \in \mathbb{N}} G_i$, where for every natural number $i$, $G_i$ is a cyclic group and the numbers $|G_i|$ are coprime, then by choosing $\sigma = (1)_{i \in \mathbb{N}}$ and $h_n := h^{\prod_{i \geq n} G_i}$, Theorem 2 applies. In fact, this is the only case: a profinite group $G$ topologically generated by an element such that for every prime number $p$, the supernatural number $p^\infty$ does not divide $\#G$ is always of this form (see \cite{Ser}, section 1.3, exercises).
\subsection{Regularity of the towers}
To prove Theorem 1, we need to understand for every $n\in\mathbb{N}$ the structure of the subfields of $H$ generated by $h_n(t,\sigma)$ and $h$. Let $t_n$ be the class of $t\in h_n(t,\sigma)$ in $H$. In particular, if $j \geq m$, then $t_j = t_m^{\frac{l_j}{l_m}}$. We will show that $h_n(t,\sigma)$ and $h$ generate a subfield of $H$ that is $h$-isomorphic to a skew fraction field with $h$ as the coefficient field.
\newline
\begin{prop}
Let $n \in \mathbb{N}$, and define $\Omega_n := h \cdot h_n(t,\sigma)$ as the compositum in $H$ of $h$ and $h_n(t,\sigma)$. We denote by $i_n$, $j_n$, and $\alpha$ the inclusion morphisms of $\Omega_n$ into $H$, from $h$ into $\Omega_n$, and from $h$ into $H$ respectively.
\begin{itemize}
\item There exists a finite order automorphism $\psi_n$ of $h$ and an isomorphism $f_n : h(t,\psi_n) \rightarrow \Omega_n$ such that $i_n \circ j_n = \alpha$.
\item The morphism $f_n \otimes 1 : h(t,\psi_n) \otimes_{\mathbf{k}(t^{l_n})} L_{l_n} \rightarrow \Omega_n \otimes_{\mathbf{k}(t^{l_n})} L_{l_n} \subset H_G$ induces an isomorphism between $h(t,\psi_n) \otimes_{\mathbf{k}(t^{l_n})} L_{l_n}$ and the compositum of $h$ and $h_n(t,\sigma) \otimes_{\mathbf{k}(t^{l_n})} L_{l_n}$ in $H_G$, and $(i_n \otimes 1) \circ (j_n \otimes 1) = \alpha \otimes 1$.
\item $h(t,\psi_n) \otimes_{\mathbf{k}(t^{l_n})} L_{l_n} / h(t,\psi_n)$ is an outer Galois extension with Galois group isomorphic to $\mathrm{Gal}(L/\mathbf{k}(t))$.
\item $H / h$ and $H_G / h$ are regular extensions.
\item We have an isomorphism of $h$-extensions $H \simeq \varinjlim h(t,\psi_n)$.
\end{itemize}
\end{prop}
\textbf{Proof:}\footnote{In this proof we use without stating it the fact that a (noncommutative)  integral domain $R$ satisfying the Ore condition have a unique (skew-)field of fractions $\mathbf{Frac}(R)$ and that any injective morphism of rings from such a ring to a field $K$ extends uniquely to a morphism of fields $\mathbf{Frac}(R)\rightarrow K$} Let $m \geq n$. If $x \in h_m$, then $t_n x = t_m^{\frac{l_m}{l_n}} x = \sigma^{\frac{l_m}{l_n}}(x) t_m^{\frac{l_m}{l_n}}$. Now, for every integer $u \in \mathbb{N}$, we have $\frac{l_{u+1}}{l_u} \equiv 1 \ (\text{mod}\ d_u)$, and for any integer $p \geq n$, we have $d_n | d_p$. Thus, if $n \leq u \leq m$, $\sigma^{\frac{l_m}{l_n}}|_{h_u} = \sigma^{\frac{l_u}{l_n}}|_{h_u}$. In particular, there exists a (unique) automorphism $\psi_n$ of $h$ whose restriction to $h_m$ is $\sigma^{\frac{l_m}{l_n}}|_{h_m}$ for all $m \geq n$. We then define a morphism 
\begin{align*}
   f_n : &h[t,\psi_n] \rightarrow H \\
   &\sum_{i=0}^k a_i t^i \mapsto \sum_{i=0}^k a_i t_n^i .
\end{align*}
The choice of $\psi_n$ ensures that this is not just a morphism of $h$-vector spaces but a ring morphism. Since the set of powers of $t_n$ is $h$-linearly independent, this is an injective morphism from $h[t,\psi_n]$ into a field, and it extends to an isomorphism between $h(t,\psi_n)$ and a subfield of $H$. This subfield of $H$ contains $h_n[t,\sigma]$ (hence $h_n(t,\sigma)$, which is its fraction field) and $h$. Conversely, any subfield of $H$ containing $h_n(t,\sigma)$ and $h$ contains the powers of $t_n$, hence the image of $f_n$. We have thus shown that the image of $f_n$ is $\Omega_n$. By construction, $i_n \circ j_n = h$.
\newline
Since $d_m | l_m$, for every $m \geq n$, we have $\sigma^{l_m}|_{h_m} = \mathrm{Id}_{h_m}$. Now, $\psi^{l_n}|_{h_m} = \sigma^{l_m}|_{h_m}$, so $\psi_n^{l_n} = \mathrm{Id}_h$, and $t^{l_n}$ is central in $h(t,\psi_n)$. Thus, $h(t,\psi_n) \otimes_{\mathbf{k}(t^{l_n})} L_{l_n}$ naturally has a structure of $\mathbf{k}(t^{l_n})$-algebra, and it is even an outer Galois extension of $h(t,\psi_n)$ with Galois group isomorphic to $\mathrm{Gal}(L/\mathbf{k}(t))$ (by Lemma 2.1.1 of \cite{BEH}). The fraction field of the $\mathbf{k}[t^{l_n}]$-algebra $h[t,\psi_n] \otimes_{\mathbf{k}[t^{l_n}]} L_{l_n}$ in $h(t,\psi_n) \otimes_{\mathbf{k}(t^{l_n})} L_{l_n}$ is $h(t,\psi) \otimes_{\mathbf{k}(t^{l_n})} L_{l_n}$: indeed, if $x \in h(t,\psi_n) \otimes_{\mathbf{k}(t^{l_n})} L_{l_n}$, there exists an $r \in \mathbb{N}$, polynomials $P_i(t), Q_i(t) \in h(t,\psi_n)$, and $R_i(t) \in L_{l_n}$ for $1 \leq i \leq r$ such that 
$$x = \sum_{i=1}^r P_i(t)^{-1} Q_i(t) \otimes R_i(t).$$
In particular, 
$$x = \sum_{i=1}^r (P_i(t)^{-1} \otimes 1)(Q_i(t) \otimes R_i(t)).$$
Thus, $P_i(t)^{-1} \otimes 1 \in (h[t,\psi_n] \otimes_{\mathbf{k}[t^{l_n}]} L_{l_n})^{-1}$ and $(Q_i(t) \otimes R_i(t)) \in h[t,\psi_n] \otimes_{\mathbf{k}[t^{l_n}]} L_{l_n}$, so $x$ belongs to the fraction field of $h[t,\psi_n] \otimes_{\mathbf{k}[t^{l_n}]} L_{l_n}$ in $h(t,\psi_n) \otimes_{\mathbf{k}(t^{l_n})} L_{l_n}$. In particular, it is enough to define an injective morphism of $h$-rings $h[t,\psi_n] \otimes_{\mathbf{k}[t^{l_n}]} L_{l_n} \rightarrow H_G$ for it to extend to a morphism of field extensions of $h$: $h(t,\psi_n) \otimes_{\mathbf{k}(t^{l_n})} L_{l_n} \rightarrow H_G$. Furthermore, since $$h[t,\psi_n] \otimes_{\mathbf{k}[t^{l_n}]} L_{l_n} = \varinjlim h_n[t,\sigma^{\frac{l_m}{l_n}}] \otimes_{\mathbf{k}[t^{l_n}]} L_{l_n}$$, it is sufficient to define a family of morphisms 
\begin{align*}f_m : h_m[t, \sigma^{\frac{l_m}{l_n}}] \otimes_{\mathbf{k}[t^{l_n}]} L_{l_n} &\rightarrow H_G \\ P(t) \otimes R(t^{l_n}) &\mapsto P(t_m^{\frac{l_m}{l_n}}) \otimes R(t^{l_m}) \in h_m[t, \sigma] \otimes_{\mathbf{k}[t^{l_m}]} L_{l_m} \subset H_G.\end{align*}
Thus, we have defined a morphism of field extensions of $h$, going from $h(t,\psi_n) \otimes_{\mathbf{k}(t^{l_n})} L_{l_n}$ into $H_G$. Its image contains $h_n(t,\sigma) \otimes_{\mathbf{k}(t^{l_n})} L_{l_n}$ and $h$, hence also the subfield generated by these. Furthermore, any subfield of $H_G$ containing $h$ and $h_n[t,\sigma] \otimes_{\mathbf{k}[t^{l_n}]} L_{l_n}$ contains the $P(t_n) \otimes 1$ for $P$ with coefficients in $h$, and the $1 \otimes R(t^{l_m}) \in 1 \otimes L_{l_m} \subset H_G$ for $R(t) \in L$. Therefore it contains the image of $f_n$. We have thus shown that the subfield of $H_G$ generated by $h_n(t,\sigma) \otimes_{\mathbf{k}(t^{l_n})} L_{l_n}$ and $h$ is $\Omega_n \otimes_{\mathbf{k}(t^{l_n})} L_{l_n}$ and is isomorphic to $h(t,\psi_n) \otimes_{\mathbf{k}(t^{l_n})} L_{l_n}$.
\newline
\newline
Let $x$ be an algebraic element of $H$ over $h$. Let $n \in \mathbb{N}$ such that $x \in h_n(t,\sigma)$. Then, by Proposition 5, there exists an automorphism $\psi$ of $h$ such that $x$ is in a subfield of $H$ isomorphic to the extension of $h$ to $h(t,\psi)$ (this field being the compositum in $H$ of $h$ and $h_n(t,\sigma)$). Now, $h(t,\psi)/h$ is a regular extension, so $x \in h$.
\newline
\newline
Let $x$ be an algebraic element of $H_G$ over $h$. Let $n \in \mathbb{N}$ such that $x \in h_n(t,\sigma) \otimes_{\mathbf{k}(t^{l_n})} L_{l_n}$. Then, by Proposition 5, there exists a finite order automorphism $\psi$ of $h$ such that $x$ is in a subfield of $H_G$ isomorphic to $h(t,\psi) \otimes_{\mathbf{k}(t^{l_n})} L_{l_n}$. Now, since $\psi$ is of finite order, there exists an $h$-embedding $h(t,\psi) \otimes_{\mathbf{k}(t^{l_n})} L_{l_n} \rightarrow h((t,\psi))$, so $h(t,\psi) \otimes_{\mathbf{k}(t^{l_n})} L_{l_n} / h$ is a regular extension and $x \in h$.
\newline
For the last point, it is simply a matter of noting that since $H$ is the increasing union of the $h_n(t,\sigma)$, it is also the increasing union of the compositums in $H$ of $h$ and the $h_n(t,\sigma)$.$\blacksquare$
\newline
\textbf{Remark:} Since $\psi_0 = \mathrm{Id}_h$, $H$ is an extension of the fraction field in a central indeterminate $h(t)$. $H / h(t)$ is even an algebraic extension of $h(t)$ since for every integer $n$, the family $(t_n^i)_{0 \leq i \leq l_n - 1}$ is a basis of the $h(t)$-vector space (both left and right) of $h(t,\psi_n) \subset H$. In particular it implies that $H$ cannot be isomorphic to $h(t,\sigma)$ if $h$ is commutative (see paragraph 2.4).
\newline
\subsection{Proof of Theorem 1}
In this subsection, we show that the generalized Artin's lemma allows us to find a subfield $k$ of $h$ and subextensions $h_n/k$ for $n \in \mathbb{N}$ satisfying the assumptions of Theorem 6 whenever the group $\mathrm{Aut}(h)$ contains a subgroup $\Gamma$ isomorphic to a product of cyclic groups of coprime orders and such that the stabilizer of any element of $h$ under the action of $\Gamma$ is open for the product topology. This point and the previous section allow us to prove Theorem 1. We then show that with the notations of Theorem 5, if the automorphism of $h$ whose restriction to $h_n$ is $\sigma^{a_n}$ generates topologically an open subgroup of $\mathrm{Gal}(h/k)$, then the result of Theorem 5 can be obtained by choosing $a_n = 1$ provided another filtration $(h_n')_{n \in \mathbb{N}}$ of $h$ is chosen.We prove a statement equivalent to Theorem 1 but with a focus on the arithmetic of the field $h$ rather than on the field extension $h/k$.
\begin{thm}
Let $h$ be a field whose automorphism group contains an infinite product $\Gamma$ of cyclic groups of coprime orders. If $h/h^\Gamma$ is an outer extension such that for every $x \in h$ the stabilizer of $x$ under the action of $\Gamma$ is an open subgroup for the product topology on $\Gamma$, $[Z(h):Z(h^\Gamma)] = +\infty$, and $Z(h^\Gamma)$ contains an ample field, then the class $\mathrm{InvReg}(h)$ contains extensions of $h$ that are not of the form $h(t,\sigma)/h$ with $\sigma$ an automorphism of finite order.
\end{thm}
\textbf{Proof}: Let $G$ be a finite group, and let $h$ be a field satisfying the hypotheses of Theorem 7. Then, according to the generalized Artin's Lemma, $h/h^\Gamma$ is an outer Galois extension with Galois group $\Gamma$. Let $k = h^\Gamma$ and $\sigma$ be a topological generator in the proof of Theorem 5. Write $\Gamma = \prod_{i \in \mathbb{N}} \mathbb{Z}/n_i \mathbb{Z}$. Let $\sigma$ be the automorphism of $h$ corresponding to $(1)_{i \in \mathbb{N}}$ and $h_i = h^{\prod_{j \geq i} \mathbb{Z}/n_i \mathbb{Z}}$ and $\sigma_i = \sigma|_{h_i}$, which satisfy the hypotheses of Theorem 6. According to Proposition 4, the extensions $H_G/h$ and $H/h$ are regular, and from the proof of Theorem 5, the extension $H_G/H$ is an outer Galois extension with a Galois group isomorphic to $G$. Therefore, $H$ belongs to $\mathrm{InvReg}(h)$.

Suppose there exists an isomorphism $f: H \to h(t, u)$ for some $u$ an automorphism of $h$ of finite order, such that the composition of $f$ and the embedding $h \subset H$ is the embedding $h \subset h(t, u)$. Since $u$ induces an automorphism of $Z(h)$, we can consider the extension $Z(h)/Z(h)^u$, and Artin’s Lemma assures us that this extension is Galois of finite degree. Now, $Z(h)^u = h \cap Z(h(t, u))$, and since $f$ induces an isomorphism between the centers and fixes $h$, we must have $[Z(h): h \cap Z(H)] < +\infty$. On the other hand, if $x \in h \cap Z(H)$, then $x$ is fixed by $\Gamma$, so $x \in h^\Gamma$. Since $x \in Z(H)$, we have that $x$ commutes with $h^\Gamma$, and thus $x \in Z(h^\Gamma)$. In particular, $h \cap Z(H) \subset Z(h^\Gamma)$, so $[Z(h): h \cap Z(H)] > [Z(h): Z(h^\Gamma)]$. Since $[Z(h): Z(h^\Gamma)] = +\infty$, we obtain a contradiction: $H$ cannot be of the form $h(t, u)$ with $u$ of finite order, which concludes the proof. Theorem 1 follows immediately because $h/k=h^\Gamma$ is an algebraic outer Galois extension of group $\Gamma$. $\blacksquare$

\textbf{Remark}: The additional condition imposed in this theorem, $[Z(h): Z(h^\Gamma)] = +\infty$, can indeed be satisfied: for example, this holds if $h$ is a commutative field.

Theorem 6 and Corollaries 1 and 2 are obtained by considering the specific case where we can choose the constant sequence $(a_n)_{n \in \mathbb{N}} = 1$ for the sequence in Theorem 5. Note that it is easy to satisfy the divisibility conditions of Theorem 6 for a filtration $h = \varinjlim h_n$ by setting $l_n := \prod_{i=0}^n d_i$ and $a_n = d_n$. However, this reduces to working at each level in a central fraction field with an indeterminate. Conversely, one can question the extent to which our choice of a constant sequence $(a_n)_{n \in \mathbb{N}} = 1$ is restrictive. In fact, we have the following result:

\begin{prop} We retain the notation of Theorem 5. Define the field $k' := \varinjlim h_n^{\sigma_n^{a_n}}$. Let $\tau$ be the unique element of $\mathrm{Gal}(h/k)$ such that for all $n \in \mathbb{N}$, we have $\tau|_{h_n} = \sigma_n^{a_n}$. If $k'/k$ is a finite-degree extension (equivalently, if $\overline{<\tau>}$ is an open subgroup of $\mathrm{Gal}(h/k)$), then for every finite group $G$, the field $H_G$ is obtained by applying Theorem 3 to $h/k'$ with $a_n = 1$ for all $n \in \mathbb{N}$. \end{prop}

\textbf{Proof}: If $k'/k$ is a finite-degree extension, there exists an integer $N$ such that $k' \subset h_N$. We have $h = \varinjlim h_{n+N}$, $\mathrm{Gal}(h/k') = \overline{<\tau>}$, and for $n \in \mathbb{N}$, we have that $\mathrm{Gal}(h_{n+N}/k')$ is cyclic and generated by $\sigma_{n+N}^{a_{n+N}}$. We have $H = \varinjlim h_{n+N}(t, \tau|_{h_{n+N}})$, $H_G = \varinjlim h_{n+N}(t, \tau|_{h_{n+N}}) \otimes_{\mathbf{k}(t^{l_{n+N}})} L_{l_{n+N}}$, $\mathbf{k} \subset Z(k')$, and the hypotheses of Theorem 5 are satisfied. Therefore, $H$ and $H_G$ are the fields obtained by replacing $k$ with $k'$, $h_n$ with $h_{n+N}$, $(d_n)_{n \in \mathbb{N}}$ with $([h_{n+N}: k'])_{n \in \mathbb{N}}$, $(a_n)_{n \in \mathbb{N}}$ with $(1)_{n \in \mathbb{N}}$, and $(l_n)_{n \in \mathbb{N}}$ with $(l_{n+N})_{n \in \mathbb{N}}$ in the statement of Theorem 5. $\blacksquare$
\printbibliography
\end{document}